\documentclass[10pt]{amsart}
\usepackage{amsmath,amscd,amssymb,epsfig}

\newcommand{\ZZ}{\mathbb{Z}}
\newcommand{\QQ}{\mathbb{Q}}

\newcommand{\Sup}{{\rm{Supp}}}
\newcommand{\Sip}{A}
\newcommand{\sign}{{\rm{sign}}}
\newcommand{\Tr}{{\rm{Tr}}}
\newcommand{\SYT}{{\rm{SYT}}}
\newcommand{\St}{{\rm{St}}}

\newcommand{\inv}{\operatorname{inv}}
\newcommand{\Inv}{\operatorname{Inv}}
\newcommand{\Pair}{{\operatorname{Pair}}}

\newcommand{\Des}{{\rm{Des}}}
\newcommand{\des}{{\rm{des}}}

\newtheorem{theorem}{Theorem}[section]
\newtheorem{corollary}[theorem]{Corollary}
\newtheorem{proposition}[theorem]{Proposition}

\newtheorem{lemma}[theorem]{Lemma}
\newtheorem{defn}[theorem]{Definition}

\newtheorem{question}[theorem]{Question}

\newtheorem{observation}[theorem]{Observation}

\numberwithin{figure}{section}

\begin{document}
\title{Combinatorial Gelfand Models}
\bibliographystyle{acm}
\author{Ron M. Adin}
\address{Department of Mathematics,
Bar-Ilan University, Ramat-Gan 52900, Israel}
\email{radin@math.biu.ac.il}
\author{Alexander Postnikov}
\address{Department of Applied Mathematics,
Massachusetts Institute of Technology, MA 02139, USA}
\email{apost@math.mit.edu}
\author{ Yuval Roichman}
\address{Department of Mathematics,
Bar-Ilan University, Ramat-Gan 52900, Israel}
\email{yuvalr@math.biu.ac.il}

\keywords{symmetric group, Iwahori-Hecke algebra, descents,
inversions, character formulas, Gelfand model}

\thanks{%First author supported by ....??
First and third authors supported in part by the Israel Science
Foundation, founded by the Israel Academy of Sciences and Humanities.}

\begin{abstract}
A combinatorial construction of a Gelfand model for the symmetric
group and its Iwahori-Hecke algebra is presented.
\end{abstract}

\date{Submitted: Oct 28, '07; Revised: March 24, '08}

\maketitle

%------------------------------------------------------------------------------
\section{Introduction}%\label{intro-section}
%------------------------------------------------------------------------------

%Denote the set of the irreducible characters of a group $G$ by
%$G^*$.
A complex representation of a group or an algebra $A$ is called a
{\it Gelfand model} for $A$, or simply a {\it model},  if it is
equivalent to the multiplicity free direct sum of {\bf all}
$A$-irreducible representations.

Models (for compact Lie groups) were first constructed by
Bernstein, Gelfand and Gelfand~\cite{BGG}. Constructions of models
for the symmetric group, using induced representations from
centralizers, were found by Klyachko~\cite{K1, K2} and by Inglis,
Richardson and Saxl~\cite{Saxl}; see also~\cite{B, R, A, A1, A3}.
Our goal is to determine an explicit and simple combinatorial
action which gives a model for the symmetric group and its
Iwahori-Hecke algebra.

%\smallskip

%------------------------------------------------------------------------------
\subsection{Signed Conjugation.}
%------------------------------------------------------------------------------

Let $S_n$ be the symmetric group on $n$ letters,
$S = \{s_1,\ldots,s_{n-1}\}$ its set of simple reflections,
$I_n = \{\pi\in S_n\,|\,\pi^2 = id\}$ its set of involutions,
and $V_n:={\rm{span}}_\QQ\{C_w\ |\ w\in I_n\}$ a vector
space over $\QQ$ formally spanned by the involutions.

\smallskip

Recall the standard length function on the symmetric group
%$\ell(\cdot)$
$$
\ell(\pi) := \min\{\ell\,|\, \pi = s_{i_1} s_{i_2} \cdots s_{i_\ell},\
s_{i_j} \in S (\forall j)\},
$$
%where $S$ is the Coxeter generating set,
the descent set
$$
\Des(\pi):=\{ s\in S\,|\, \ell(\pi s) <\ell(\pi)\},
$$
and the descent number ${\rm{des}}(\pi):=\#\Des(\pi)$.

\medskip

Define a map $\rho:S\to GL(V_n)$ by
\begin{equation}\label{e.signed_conjugation2}
\rho(s)C_w := {\rm{sign}}(s;\ w)\cdot C_{sws}\qquad (\forall s\in
S, w\in I_n)
\end{equation}
where
\begin{equation}\label{e.sign-definition}
{\rm{sign}}(s;\ w):=
\begin{cases}
   -1, & \text{ if } sws=w \text { and } s\in \Des(w);\\
   1, & \text{ otherwise. }
\end{cases}
\end{equation}

\begin{theorem}\label{t.representation}
$\rho$ determines an $S_n$-representation. \end{theorem}

\begin{theorem}\label{t.model}
$\rho$ determines a Gelfand model for  $S_n$.
%; namely, it is equivalent as an $S_n$ module to
%the multiplicity free sum of all $S_n$ irreducible
%representations.
\end{theorem}

%\medskip

%------------------------------------------------------------------------------
\subsection{Hecke Algebra Action.}%\label{Intro-Hecke}
%------------------------------------------------------------------------------

Consider ${H}_n(q)$, the Hecke algebra of the symmetric group $S_n$
(say over the field $\QQ(q^{1/2})$),
with set of generators $\{T_i\,|\, 1\le i<n\}$ and defining relations
$$
(T_i+q)(T_i-1)=0\qquad (\forall i),
$$
$$
T_iT_j=T_jT_i\qquad {\rm if\ \ } |i-j|>1,
$$
$$
T_iT_{i+1}T_i=T_{i+1}T_iT_{i+1}\qquad (1\le i< n-1) .
$$
Note that some authors use a slightly different notation, with $T_i$
consistently replaced by $-T_i$.

\smallskip

In order to construct an extended signed conjugation providing %, which gives
a model for ${H}_n(q)$, we extend the standard notions of length
and weak order. Recall that the (left) weak order on $S_n$ is the
reflexive and transitive closure of the relation:
$w \prec_L ws$ if $s \in S$ and $\ell(sw) = \ell(w)+1$.

\begin{defn}\label{involutive}%\rm\\
Define the {\em involutive length} of an involution $w\in I_n$
of cycle type $2^k 1^{n-2k}$ as
$$
\hat\ell(w):=\min\{\ell(v)|\ w=vs_1s_3 \cdots s_{2k-1}v^{-1},\,v\in S_n\},
$$
where $\ell(v)$ is the standard length of $v\in S_n$.

Define the {\em involutive weak order} $\le_I$ on $I_n$ as the
reflexive and transitive closure of the relation: $w \prec_I sws$ %\lessdot
if $s \in S$ and $\hat\ell(sws) = \hat\ell(w)+1$.
\end{defn}

\bigskip

Now define a map $\rho_q:S\to GL(V_n)$ by
\begin{equation}\label{e.Hecke-action}
\rho_q(T_s)C_w :=  \begin{cases}
   -q C_w,   & \text{ if } sws=w \text { and } s\in \Des(w); \\
   C_w,      & \text{ if } sws=w \text{ and } s\not\in \Des(w); \\
   (1-q) C_w + q C_{sws}, & \text{ if } w <_I sws; \\
   C_{sws},  & \text{ if } sws <_I w.
\end{cases}
\end{equation}
%where $\Des(\cdot)$ is the standard descent set and $<_I$ is
%the involutive weak order. % defined above.

%\begin{theorem}\label{t.Hecke-representation} $\rho_q$ is an
%${H}_n(q)$-representation. \end{theorem}

\begin{theorem}\label{t.Hecke-model}
$\rho_q$ is a Gelfand model for $H_n(q)$ ($q$ indeterminate); namely,
\begin{itemize}
\item[(1)]
$\rho_q$ is an ${H}_n(q)$-representation.
\item[(2)]
$\rho_q$ is  equivalent %as an ${H}_n(q)$-module
 to the multiplicity free sum of all irreducible ${H}_n(q)$-representations.
\end{itemize}

\end{theorem}

The proof involves Lusztig's version of Tits' deformation theorem~\cite{L}.
For other versions of this theorem see~\cite[\S 4]{Bou},
\cite[\S 68.A]{CR} and~\cite{BC}.
%, recalling that, by Theorem~\ref{t.model},
%$\rho_{|q=1}=\rho$ is a model for the group algebra of $S_n$.

\medskip

%\begin{remark}\rm
Let $\mu=(\mu_1,\mu_2,\dots,\mu_t)$ be a partition of $n$ and let
$a_j:=\sum_{i=1}^j \mu_i$ $(1\le j\le t)$. A permutation $\pi\in S_n$ is
{\it $\mu$-unimodal} if for every $0\le j< t$ there exists
$1\le d_j\le \mu_{j+1}$ such that
$$
\pi_{a_j+1} < \pi_{a_j+2} < \ldots < \pi_{a_j+d_j}
> \pi_{a_j+d_j+1}> \ldots > \pi_{a_{j+1}}.
$$
The character of $\rho_q$ may be expressed as a generating function
for the descent number over $\mu$-unimodal involutions.

\begin{proposition}\label{t.Hecke-character}
$$
\Tr(\rho_q(T_\mu))= \sum\limits_{\{w\in I_n\,|\, w\text{ is $\mu$-unimodal}\}}
(-q)^{\des(w)}
$$
where
$$
T_\mu:=T_1 T_2 \cdots T_{\mu_1-1} T_{\mu_1+1} \cdots T_{\mu_1+\ldots+\mu_t-1}
$$
is the subproduct of $T_1 T_2 \cdots T_{n-1}$ obtained by
omitting $T_{\mu_1+\dots+\mu_i}$ for all $1\le i<t$.
\end{proposition}
%\end{remark}

%------------------------------------------------------------------------------
\section{Proof of Theorem~\ref{t.representation}}
%\section{The Involutive Inversion Number}
%\label{sections-inversions}
%------------------------------------------------------------------------------

%------------------------------------------------------------------------------
\subsection{First Proof}
%------------------------------------------------------------------------------

This proof relies on a variant of the inversion number, which is
introduced in this section. Recall the definition of the inversion
set of a permutation $\pi\in S_n$,
$$
\Inv(\pi):=\{\,\{i,j\} \,|\, (j-i) \cdot (\pi(j)-\pi(i)) < 0 \}.
$$

%, yielding another presentation of
%$S_n$-model as a signed conjugation.
%\medskip

\begin{defn}%\label{Pair-Inv}
For an involution $w\in I_n$ let $\Pair(w)$ be the set of 2-cycles of $w$
(considered as unordered 2-sets).
For a permutation $\pi\in S_n$ and an involution $w\in I_n$ let
$$
\Inv_w(\pi):= \Inv(\pi) \cap \Pair (w)
$$
and
$$
\inv_w(\pi):=\# \Inv_w(\pi).
$$
\end{defn}

Now redefine $\rho:S_n\to GL(V_n)$ by
%\medskip
\begin{equation}\label{e.alternative-signed-conjugation}
\rho(\pi)C_w := (-1)^{\inv_w(\pi)}\cdot C_{\pi w\pi^{-1}}\qquad
(\forall \pi\in S_n, w\in I_n).
\end{equation}
%\medskip
%
%\begin{proof} {\bf - of Theorem~\ref{t.representation}.}
Note that for every Coxeter generator $s = (i,i+1)\in S$ and every
involution $w\in I_n$,
%$\inv_w(s)=1$ if and only if   $s\in
%\Des(w)$ and $sws=w$, and $\inv_w(s)=0$ otherwise.
\begin{eqnarray*}
\inv_w(s) &=&
\begin{cases}
   1,   & \text{ if } w(i) = i+1; \\
   0,   & \text{ otherwise}
\end{cases}\\
&=&
\begin{cases}
   1,   & \text{ if } sws=w \text { and } s\in \Des(w); \\
   0,   & \text{ otherwise. }
\end{cases}
\end{eqnarray*}
Thus, definition~(\ref{e.alternative-signed-conjugation}) of $\rho$
coincides on the Coxeter generators with
the original definition~(\ref{e.signed_conjugation2}).
In order to prove that $\rho$ is an $S_n$-representation
%which coincide on all elements in $S_n$,
it suffices to prove that $\rho$ is a group homomorphism.

Indeed, for every pair of permutations $\sigma,\pi\in S_n$, every
involution $w\in I_n$, and every $1\le i< j\le n$,
%$(i,j)\in \Inv_w(\sigma\pi)$ if and only if
%$(i,j)\in \Pair (w)$ and
%%${\sigma\pi(j)-\sigma\pi(i)\over
%%j-i}<0$,
%%if and only if $(\pi(i),\pi(j))\in \Pair(\pi w\pi^{-1})$
%$$ 0>{\sigma\pi(j)-\sigma\pi(i)\over j-i}={\pi(j)-\pi(i)\over
%j-i}\cdot {\sigma(\pi(j))-\sigma(\pi(i))\over \pi(j)-\pi(i)}.$$
%Thus
$$
\chi[\{i,j\}\in \Inv_w(\sigma\pi)] =
\chi[\{i,j\}\in\Inv_w(\pi)] \cdot
\chi[\{\pi(i),\pi(j)\}\in\Inv_{\pi w\pi^{-1}}(\sigma)],
$$
where $\chi[\text{ event }] := -1$ if the event holds and $1$ otherwise.
Hence, for every pair of
permutations $\sigma,\pi\in S_n$ and every involution $w\in I_n$,
$$
(-1)^{\inv_w(\sigma\pi)}=(-1)^{\inv_w(\pi)}\cdot(-1)^{\inv_{\pi
w\pi^{-1}}(\sigma)},
$$
and thus
\begin{eqnarray*}
\rho(\sigma\pi)C_w
&=& (-1)^{\inv_w(\sigma\pi)}\cdot C_{(\sigma\pi) w(\sigma\pi)^{-1}}\\
&=& (-1)^{\inv_w(\pi)} \cdot (-1)^{\inv_{\pi w\pi^{-1}}(\sigma)}C_{\sigma(\pi w\pi^{-1})\sigma^{-1}} \\
&=& (-1)^{\inv_w(\pi)} \cdot \rho(\sigma)(C_{\pi w\pi^{-1}}) = \rho(\sigma)(\rho(\pi)C_w).
\end{eqnarray*}
This proves that $\rho$ is an $S_n$-representation,
completing the proof of Theorem~\ref{t.representation}.

\qed

%\end{proof}

%We summarize

%\begin{corollary}\label{inversions-theorem}
%$\rho$ and $\psi$ are identical $S_n$-representations.
%\end{corollary}

%------------------------------------------------------------------------------
\subsection{Second Proof.}% of Theorem~\ref{t.representation}}
%------------------------------------------------------------------------------

In order to prove that $\rho$ (defined on $S$) extends to an $S_n$-representation it
suffices to verify the relations:
%\begin{equation}\label{relation1}
$$
\rho(s)^2=1\qquad (\forall s\in S),
$$
%\end{equation}
%\begin{equation}\label{relation2}
$$
\rho(s)\rho(t)=\rho(t)\rho(s)\qquad {\rm if\ \ } st=ts,
$$
%\end{equation}
%\begin{equation}\label{relation3}
$$
\rho(s)\rho(t)\rho(s)=\rho(t)\rho(s)\rho(t)\qquad {\rm if\ \ } sts=tst.
$$
%\end{equation}

We will prove the third relation.
%Verification of the first two relations is easy and left to the reader.
Verifying the other two relations is easier and will be left to the reader.

Let $s=(i,i+1)$ and $t=(i+1,i+2)$. For every permutation $\pi\in S_n$ let
$$
\Sup(\pi) := \{i\in [n]\ |\ \pi(i)\ne i\}.
$$
Denote by $O(w)$ the orbit of an involution $w$
under the conjugation action of $\langle s,t\rangle$,
the subgroup of $S_n$ generated by $s$ and $t$.
Since $w$ is an involution $\#O(w)\ne 2$; hence there are three options $\#O(w)=1,3,6$.

\bigskip

\noindent{\bf Case (a).} $\#O(w)=1$. Then $sws=w$ and $twt=w$.
Furthermore, in this case $\Sup(w)\cap \{i,i+1,i+2\}=\emptyset$,
so that $\sign(s;\ w)=\sign(t;\ w)=1$; thus
$\rho(s)\rho(t)\rho(s)C_w=\rho(t)\rho(s)\rho(t)C_w=C_w$.

\bigskip

\noindent{\bf Case (b).} $\#O(w)=3$. (This happens, for example, when $w=s$.)
With no loss of generality there exists an element $v$ in the orbit such that
%\begin{equation}\label{e.condition1}
$$
v, tvt, stvts\ \ {\rm \ are\ distinct\ elements\ in\ the\ orbit,}
$$
%\end{equation}
while
\begin{equation}\label{e.condition2}
svs=v \ \ {\rm and }\ \ t(stvts)t=stvts .
\end{equation}
Thus
$$
\rho(s)= \left(\begin{matrix}
x         & 0         & 0       \\
0                  & 0         & 1   \\
0                  & 1         & 0   \\
\end{matrix}\right)
$$
and
$$
\rho(t)=\left(
\begin{matrix}
0        & 1         & 0       \\
1                  & 0         & 0   \\
0                  & 0         &  z  \\
\end{matrix}\right),
$$
where $x=\sign(s;\ v)$ and $z=\sign(t;\ stvts)$.
$\rho(s)\rho(t)\rho(s)=\rho(t)\rho(s)\rho(t)$ holds if and only if
$x=z$, which holds if and only if
\begin{equation}\label{e.iff8}
s\in \Des(v)\Longleftrightarrow t\in\Des(stvts).
\end{equation}
To prove this, observe that for every $w\in S_n$ and $s\in S$ the
following holds :
\begin{itemize}
\item[$(A)$] $sws=w$  and $s\not\in\Des(w)$ if and only if
$\Sup(w)\cap\Sup(s)=\emptyset$. \item[$(B)$] $sws=w$  and
$s\in\Des(w)$ if and only if $w=us$, where
$\Sup(u)\cap\Sup(s)=\emptyset$.
\end{itemize}
Assuming $t\not\in\Des(stvts)$ implies, by (\ref{e.condition2}) and
$(A)$, that $\Sup(stvts)\cap \Sup(t)=\emptyset$. Hence
$$
stvts(i+1)=i+1.
$$
On the other hand, assuming $s\in \Des(v)$ implies, by
(\ref{e.condition2}) and  $(B)$, that there exists $u=vs$ with
$i+1\not\in \Sup(u)$. Hence
$$
stvts(i+1)=stusts(i+1)=i+2,
$$
a contradiction. Similarly, assuming $s\not\in\Des(v)$ and
$t\in\Des(stvts)$ yields a contradiction (to verify this, replace
$v$ by $stvts$ and $s$ by $t$). This completes the proof of Case
(b).

\bigskip

\noindent{\bf Case (c).} $\#O(w)=6$ (this occurs, for example, when
$s=(i,i+1), t=(i+1,i+2)$ and $w=(i,j)(i+1,k)$ where $j,k\ne i+2$).
Then, for every element $v$ in the orbit, $svs\ne v$ and $tvt\ne v$.
It follows that
$$
\rho(s)\rho(t)\rho(s)C_w=C_{stswsts}=C_{tstwtst}=\rho(t)\rho(s)\rho(t)C_w.
$$
This completes the proof of the third relation .

\qed

%------------------------------------------------------------------------------
\section{Characters}
%\section{Proof of Theorem~\ref{t.model}}
%\label{Section3}
%------------------------------------------------------------------------------

%------------------------------------------------------------------------------
\subsection{Character Formula}%\label{explicit-formula}
%------------------------------------------------------------------------------

The following classical result %will be used in the proof. is due to
follows from the work of Frobenius and Schur, see~\cite[\S 4]{Is}
and ~\cite[\S 7, Ex. 69]{EC2}.

\begin{theorem}\label{t.roots}%\cite[\S 7, Ex. 69]{EC2}
Let $G$ be a finite group, for which every complex representation
is equivalent to a real representation. Then for every $w\in G$
$$
\sum\limits_{\chi\in G^*} \chi(w)=\# \{u\in G\ |\ u^2=w\},
$$
where $G^*$ denotes the set of the irreducible characters of $G$.
\end{theorem}

%Let $I_{G}:=\{w\in G\ |\ w^2=id\}$ be the set of involutions in a
%group $G$. Motivated by theorem~\ref{t.roots}, it is desirable to
%find a linear action   (of a group $G$ whose representations are
%real) on $I_G$, which realizes a model.

It is well known~\cite{Springer} that all complex representations
of a Weyl group are equivalent to rational representations. In
particular,  Theorem~\ref{t.roots} holds for $G=S_n$. One concludes

\begin{corollary}\label{formula}
Let $\pi\in S_n$ have cycle structure $1^{d_1}2^{d_2}\cdots
n^{d_n}$. Then
$$
\sum\limits_{\chi\in {S_n}^*} \chi(\pi)=
%\begin{cases}
%   0   & \text{ if there exists an even } r \text { with odd } d_r \\
   \prod\limits_{r=1}^n f(r,d_r), %& \text{ otherwise }
%\end{cases} ,
$$
where
$$
f(r,d_r) :=
\begin{cases}
   0,   & \text{ if } r \text{ is even and } d_r \text{ is odd};\\
%   1,   & \text{ if } d_r=0; \\
   {d_r\choose 2, \dots,2} \cdot r^{d_r/2},   & \text{ if } r \text{ and } d_r \text{ are even}; \\
   \sum\limits_{k=0}^{\lfloor d_r/2\rfloor} {d_r\choose d_r-2k,2,2,\dots,2}\cdot r^k, &
   \text{ if } r \text{ is odd.}
\end{cases}
$$
In particular, $f(r,0) = 1$ for all $r$.
\end{corollary}

\begin{proof}
For every $A\subseteq [n]$ let
$$
S_A:=\{\pi\in S_n\,|\,\Sup(\pi)\subseteq A\}
$$
be the subgroup of $S_n$ consisting of all the permutations whose support is contained in $A$.
%
%\smallskip
%
For every $\pi\in S_n$ and $1\le r\le n$ let $\Sip(\pi,r)\subseteq [n]$ be
the set of all letters which appear in cycles of length $r$ in $\pi$.
In other words,
$$
\Sip(\pi,r):=\{i\in [n]\,|\, \pi^r(i)=i \text{ and } (\forall j<r)\,\pi^j(i)\ne i\}
$$
For example, $\Sip(\pi,1)$ is the set of fixed points of $\pi$.

Denote by $\pi_{|r}$ the restriction of $\pi$ to $\Sip(\pi,r)$.
Then $\pi_{|r}$ may be considered as a permutation in
$S_{\Sip(\pi,r)}$.

\begin{observation}%\label{obs1}
For every $\pi\in S_n$
$$
\{u\in S_n\,|\, u^2=\pi\} =
\prod\limits_{r\ge 1} \{u_r\in S_{\Sip(\pi,r)}\,|\, u_r^2=\pi_{|r}
% \ {\rm{ and }}\ \Sup(w)\subseteq \Sip(\pi,r)
\}.
$$
\end{observation}

\begin{observation}%\label{lemma1}
Let $\pi\in S_n$ have cycle type $r^{n/r}$. Then
$$
\#\{u\in S_n\ |\ u^2=\pi\} =
\begin{cases}
   0,   & \text{ if } $r$ \text { is even and } n/r \text{ is odd}; \\
   {n/r \choose 2,\dots,2}\cdot r^{n/2r},   & \text{ if } r \text { and } n/r \text{ are even}; \\
   \sum\limits_{k=0}^{\lfloor n/2r\rfloor} {n/r\choose n/r-2k,2,2,\dots,2}\cdot r^k,
      & \text{ if } $r$ \text{ is odd }.
\end{cases}
$$
\end{observation}

Combining these observations with Theorem~\ref{t.roots} implies Corollary~\ref{formula}.

\end{proof}

\medskip

%------------------------------------------------------------------------------
\subsection{Proof of Theorem~\ref{t.model}} %\label{intro-section}
%------------------------------------------------------------------------------

%By Theorem~\ref{t.roots}, it suffices to prove that for every
%permutation $\pi\in S_n$
%$$
%\Tr (\psi(\pi))=\# \{w\in S_n\ |\ w^2=\pi\}.
%$$

%Now
We shall compute the character of the representation $\rho$ and
compare it with Corollary~\ref{formula}.
By~(\ref{e.alternative-signed-conjugation}),
%Corollary~\ref{inversions-theorem},
%\begin{equation}\label{eq0}
$$
\Tr (\rho(\pi))=\sum\limits_{w\in I_n \cap \St_n(\pi)}
(-1)^{\inv_w(\pi)},
$$
%\end{equation}
where $\St_n(\pi)$ is the stabilizer of $\pi$ under the conjugation action
of $S_n$ (i.e., the centralizer of $\pi$ in $S_n$).

\begin{observation}\label{obs.cycles}
Let $\pi\in S_n$, $w\in I_n\cap \St_n(\pi)$ and $a_1\in [n]$ any letter. %$\Sup(w)$.
Then one of the following holds:
\begin{itemize}
\item[(1)]
%$(a_1,a_2)$ is a cycle in $w$ and $a_1,a_2\not\in\Sup(\pi)$.
$(a_1,a_2,\ldots,a_r)$ is a cycle in $\pi$ $(r\ge 1)$;
$a_1$, $a_2$, $\ldots$, $a_r$ are fixed points of $w$.
\item[(2)]
$(a_1,a_2,\ldots,a_r)$ and $(a_{r+1},\ldots,a_{2r})$ are cycles in $\pi$ $(r\ge 1)$;
$(a_1,a_{r+1})$, $(a_2,a_{r+2})$, $\ldots$, $(a_r,a_{2r})$ are cycles in $w$.
\item[(3)]
$(a_1,a_2,\ldots,a_{2m})$ is a cycle in $\pi$ $(m\ge 1)$;
$(a_1,a_{m+1})$, $(a_2,a_{m+2})$, $\ldots$, $(a_m,a_{2m})$ are cycles in $w$.
\end{itemize}
\end{observation}

It follows that

%By Observation~\ref{obs.cycles},

\begin{corollary}%\label{obs.union}
%For every $\pi\in S_n$
%$$
%\St_n(\pi)\cap I_n= \prod\limits_{r=1}^n
%(C_{S_{\Sip(\pi,r)}}(\pi_{|r})\cap I_{S_{\Sip(\pi,r)}}),
%$$
%where $\Sip(\pi,r)$, $\pi_{|r}$ and $S_{\Sip(\pi,r)}$ are defined
%as in the proof of Corollary~\ref{formula}.\\
%Furthermore,
%In other words,
Fix $\pi\in S_n$. Each $w\in I_n \cap \St_n(\pi)$ has a unique decomposition
$$
%w=w_1\cdots w_r\cdots w_n
w = \prod_{r\ge 1} w_r,
$$
where
$$
w_r\in I_{S_{\Sip(\pi,r)}} \cap \St_{S_{\Sip(\pi,r)}}(\pi_{|r}) \qquad(\forall r)
$$
and $\Sip(\pi,r)$, $\pi_{|r}$ and $S_{\Sip(\pi,r)}$ are defined
as in the proof of Corollary~\ref{formula}; and
$$
\Inv_w(\pi)=\bigcup_{r\ge 1} \Inv_{w_r}(\pi_{|r}), %\biguplus
$$
a disjoint union.
%$$
%\sign(\pi;\ w)=\prod\limits_{r=1}^n \sign(\pi_{|r};\ w_r).
%$$
\end{corollary}

%Thus the theorem is reduced to proving
Hence, it suffices to prove that $\Tr(\rho(\pi))$ is equal to the
right hand side of the formula in Corollary~\ref{formula},
for $\pi$ of cycle type $r^{n/r}$.
%Without loss of generality
Since $\rho$ is a class function, we may assume that
\begin{equation}\label{e.pi}
\pi=(1,2,\dots,r)(r+1,\dots,2r)\cdots(n-r+1,n-r+2,\dots,n).
\end{equation}

%By Corollary~\ref{obs.union} and Observation~\ref{obs1}, it suffices to
%prove (\ref{eq1}) on permutations of cycle type $r^{n/r}$.

\begin{observation}\label{obs.sign}
Let $r$ be a positive integer.
\begin{itemize}
\item[(1)]
If $i$ and $j$ are distinct nonnegative integers,
%$\pi = (ir+1,ir+2,\ldots,ir+r)(jr+1,jr+2,\ldots,jr+r)$ and
$\pi$ as in (\ref{e.pi}) above, and
$w = (ir+1,jr+\sigma(1))(ir+2,jr+\sigma(2)) \cdots (ir+r,jr+\sigma(r))$
(where $\sigma$ is some power of the cyclic permutatation $(1,2,\ldots,r)$),
then
%$$
%\sign(\ (a_1,a_2,\dots,a_r)(a_{r+1},\dots,a_{2r})\ ;\
%(a_1,a_{r+1})(a_2,a_{r+2})\cdots (a_r,a_{2r})\ )=1.
%$$
%$$
%\sign(\ (1,2,\dots,r)({r+1},\dots,{2r})\ ;\
%(1,{r+1})(2,{r+2})\cdots (r,{2r})\ )=1.
%$$
$$
(-1)^{\inv_w(\pi)}=1.
$$
\item[(2)]
If $r=2m$ is even,
%$\pi = (1,2,\ldots,2m)$ and
$\pi$ as in (\ref{e.pi}) above, and
$w = (1,m+1)(2,m+2)\cdots(m,2m)$,
then
$$
(-1)^{\inv_w(\pi)}=-1.
$$
\end{itemize}
\end{observation}

\begin{lemma}\label{lemma3}
For every odd $r$ and a permutation $\pi$ as in~(\ref{e.pi}) above,
%$\pi=(1,2,\dots,r)(r+1,\dots,2r)\cdots(n-r+1,n-r+2,\dots,n)$,
%$\pi\in S_n$ be of cycle type $r^{n/r}$. Then
$$
\sum\limits_{w\in I_n \cap \St_n(\pi)} (-1)^{\inv_w(\pi)}=
\#\ (I_n \cap \St_n(\pi)) = \sum\limits_{k=0}^{\lfloor n/2r\rfloor}
{n/r\choose n/r-2k,2,2,\dots,2}\cdot r^k. % = \# I_{S_{n/r}},
$$
%the number of matchings on $n/r$ points.
\end{lemma}

\noindent{\bf Proof of Lemma~\ref{lemma3}.}
If $r$ is odd then only cases (1) and (2) in Observation~\ref{obs.cycles} are possible.
The first equality in the statement of the lemma then follows from Observation~\ref{obs.sign}(1).
The second equality follows from Observation~\ref{obs.cycles}(1)(2), counting
the involutions $w\in I_n \cap \St_n(\pi)$ with $\#\,\Sup(w) = 2rk$.

\qed

\begin{lemma}\label{lemma4}
For every even $r$ and a permutation $\pi$ as in~(\ref{e.pi}) above,
%Let $r$ be even and $\pi\in S_n$ be of cycle type $r^{n/r}$. Then
$$
%\#\ \St_n\cap I_n =
\sum\limits_{w\in I_n \cap \St_n(\pi)} (-1)^{\inv_w(\pi)} =
\begin{cases}
   0,   & \text{ if } $n/r$ \text{ is odd}; \\
   {n/r \choose 2,\dots,2}\cdot r^{n/2r},    & \text{ if } $n/r$ \text{ is even}.
   %\\ \# I_{S_{n/r}} & \text{ if } $r$ \text{ is odd }.
\end{cases}
$$
\end{lemma}

\noindent{\bf Proof of Lemma~\ref{lemma4}.}
Let $c_i = (ir+1,ir+2,\ldots,ir+r)$ be one of the cycles of $\pi$,
as in~(\ref{e.pi}). By Observation~\ref{obs.cycles},
an involution $w\in I_n\cap \St_n(\pi)$ has one of the following
three types with respect to $c_i$:
\begin{description}
\item[\bf Type (1)]
Each element of $c_i$ is a fixed point of $w$.
\item[\bf Type (2)]
$w$ maps $c_i$ onto a different cycle of $\pi$.
\item[\bf Type (3)]
$r = 2m$ is even, and $c_i$ is a union of 2-cycles of $w$:
$$
\{ir+t,ir+t+m\}\in\Pair(w)\qquad(1\le t\le m).
$$
\end{description}
Denote
$$
P_2 := \{w\in I_n\cap \St_n(\pi)\,|\, w \text{ is of type (2) w.r.t.\ all cycles of } \pi\}.
$$
For any $w\in (I_n\cap \St_n(\pi)) \setminus P_2$, let
$$
i(w) := \min\{i\,|\, w \text{ is of type (1) or (3) w.r.t.\ the cycle } c_i\}.
$$
Denote
$$
P_1 := \{w\in (I_n\cap \St_n(\pi)) \setminus P_2\,|\,
       w \text{ is of type (1) w.r.t.\ the cycle } c_{i(w)}\}
$$
and
$$
P_3 := \{w\in (I_n\cap \St_n(\pi)) \setminus P_2\,|\,
       w \text{ is of type (3) w.r.t.\ the cycle } c_{i(w)}\}.
$$
The map $\varphi: P_1 \to P_3$ which changes the action of $w$ on $c_{i(w)}$
from type (1) to type (3) is clearly a well-defined bijection;
and, by Observation~\ref{obs.sign}(2),
it reverses the sign of $(-1)^{\inv_w(\pi)}$.
The contributions of $P_1$ and $P_3$ to the sum therefore cancel each other.
Each element of the remaining set $P_2$ contributes $1$, by Observation~\ref{obs.sign}(1).
%By Observation~\ref{obs.cycles},
%the following two sets of elements in $I_n\cap \St_n(\pi)$
%%, which fix the support
%%of the $k$-th $r$-cycle of $\pi$
%$$
%A_k := \{w\in \St_n(\pi)\ : \ (k-1)r+1,(k-1)r+2,\dots,kr \not\in\Sup(w)\},
%$$
%and
%$$
%B_k := \{w\in \St_n(\pi)\ : \ \{i,i+r\} \in \Pair(w) \text{ for } (k-1)r+1\le i\le (k-1)r+r/2\}.
%$$
%% the set of involutions in
%%$\St_n(\pi)$, for which $(i,r/2+i)$, $kr+1\le i\le kr+r/2$,
%%are 2-cycles.
%Clearly, these two sets have the same cardinality.
%By Observation~\ref{obs.sign}(2), their contributions $(-1)^{\inv_w(\pi)}$
%have opposite signs.
%We are left with involutions in $\St_n(\pi)$, for which all
%2-cycles are of second type in Observation~\ref{obs.cycles}.
Lemma~\ref{lemma4} follows.

\qed

%Combining (\ref{eq0}) and Corollary~\ref{obs.union} with
%Lemmas~\ref{lemma3} and~\ref{lemma4} shows that $\Tr(\rho(\pi))$
%is equal to the right hand side of Corollary~\ref{formula},
%completing the proof of Theorem~\ref{t.model}.

Lemmas~\ref{lemma3} and~\ref{lemma4} complete the proof of
Theorem~\ref{t.model}.

\qed

%------------------------------------------------------------------------------
\section{The Hecke Algebra}%\label{Section4}
%\section{Proofs of Theorems~\ref{t.Hecke-representation}
%and~\ref{t.Hecke-model}}
%------------------------------------------------------------------------------

%------------------------------------------------------------------------------
\subsection{A Combinatorial Lemma.}
%------------------------------------------------------------------------------

Recall %the definition of the involutive length $\hat\ell$ from
Definition~\ref{involutive}. In order to prove
Theorem~\ref{t.Hecke-model} we need the following combinatorial
interpretation of the involutive length $\hat\ell$.

\begin{lemma}\label{t.ell}
Let $w\in S_n$ be an involution of cycle type $2^k1^{n-2k}$. Then
\begin{equation}\label{e.interp}
\hat\ell(w) :=
\left[\sum\limits_{t\in \Sup(w)} t - {2k+1\choose 2}\right]
+ \frac{1}{2} \left[\inv(w_{|\Sup(w)}) - k\right].
\end{equation}
\end{lemma}

%Recall Definition~\ref{Pair-Inv}.
%The following definition is used in the below proofs.
%
%\begin{defn} Let $w$ be  an involution
% of cycle type $2^k1^{n-2k}$. Then $w$ is a product of $k$ disjoint transpositions
% $t_1\cdots t_k$, and we say that each of these transpositions {\em
%factors} $w$.
%\end{defn}

\begin{proof}
Denote the right hand side of (\ref{e.interp}) by $f(w)$.
It is easy to verify that $f(w) = 0$ when $\hat\ell(w) = 0$,
i.e., when $w = s_1s_3 \cdots s_{2k-1}$.
Let $u$ and $v=s_ius_i$ be involutions in $S_n$ with
$\hat\ell(v)=\hat\ell(u)+1$. Then $|\{i,i+1\}\cap \Sup(u)|>0$.
If $|\{i,i+1\}\cap \Sup(u)|=1$ then
$$
\sum\limits_{t\in \Sup(v)} t - \sum\limits_{t\in \Sup(u)} t = \pm 1
$$
and
$\inv (v_{|\Sup(v)})=\inv (u_{|\Sup(u)})$.
If $|\{i,i+1\}\cap \Sup(u)|=2$ then
$$
\sum\limits_{t\in \Sup(v)} t = \sum\limits_{t\in \Sup(u)} t
$$
and $\inv (v_{|\Sup(v)})-\inv (u_{|\Sup(u)})\in \{2, 0, -2\}$.
Thus in both cases $|f(v)-f(u)|\le 1$.
This proves, by induction on $\hat\ell$, that $f(w)\le \hat\ell(w)$
for every involution $w$.

On the other hand, if $w$ is an involution with $f(w)>0$ then
either $\sum\limits_{t\in \Sup(w)} t>{2k+1\choose 2}$, or
$\sum\limits_{t\in \Sup(w)} t={2k+1\choose 2}$ and $\inv
(w_{|\Sup(w)})>k$. In the first case there exists $i+1\in\Sup(w)$
such that $i\not\in\Sup(w)$. Then $f(s_iws_i)=f(w)-1$. In the
second case $\Sup(w) = \{1,\ldots,2k\}$. Since $\inv
(w_{|\Sup(w)})>k$, $w\ne s_1s_3\cdots s_{2k-1}$. Thus there must
be a minimal $i$ such that $w(i)>i+1$. Let $j:= w(i)-1$ then
$w(j)>w(j+1)=i$, so $f(s_j w s_j)= f(w)-1$.
%Letting $i$ be the smallest letter such that $w(i) - i > 1$, we
%must have $w(i),w(i+1) > i+1$. Thus there exist $i+1 < j < k$ such
%that one of the following holds: either $\{i,j\}$, $\{i+1,k\}\in
%\Pair(w)$, and then $f(s_{i+1}ws_{i+1})=f(w)-1$; or $\{i,k\}$,
%$\{i+1,j\}\in\Pair(w)$, and then $f(s_iws_i)=f(w)-1$.
We conclude that $\hat\ell(w)\le f(w)$ for every involution $w$.

\end{proof}

%------------------------------------------------------------------------------
%\subsection{Proof of Theorem~\ref{t.Hecke-representation}}
\subsection{Proof of Theorem~\ref{t.Hecke-model}}
%------------------------------------------------------------------------------

The proof consists of two parts. In the first part we prove that
$\rho_q$ is an $H_n(q)$-representation by verifying the defining
relations along the lines of the second proof of
Theorem~\ref{t.representation}. In the second part we apply
Lusztig's version of Tits' deformation theorem to prove that
$\rho_q$ is a Gelfand model.

\bigskip

\noindent{\bf Part 1: Proof of Theorem~\ref{t.Hecke-model}(1).}
%$\rho_q$ is an $H_n(q)$-representation.
First, consider the braid relation
$T_iT_{i+1}T_i=T_{i+1}T_iT_{i+1}$. To verify this relation observe
that there are six possible types of orbits of an involution $w$
under conjugation by $\langle s_i,s_{i+1}\rangle$, the subgroup of $S_n$
generated by $s_i$ and $s_{i+1}$. These orbits differ by
the action of $w$ on the letters $i,i+1,i+2$ :
\begin{itemize}
\item[1.] $i,i+1,i+2\not\in \Sup(w)$. %This is an orbit of order 1.
\item[2.] Exactly one of the letters $i,i+1,i+2$ is in $\Sup(w)$.
%This is an orbit of order 3.
\item[3.] Exactly two of the letters $i,i+1,i+2$ are in $\Sup(w)$,
and these two letters form a 2-cycle in $w$.
%transposition in $\Pair(w)$.  %This is an orbit of order 3.
\item[4.] Exactly two of the letters $i,i+1,i+2$ are in $\Sup(w)$,
and these two letters do not form a 2-cycle in $w$.
%transposition in $\Pair(w)$. %This is an orbit of order 6.
\item[5.] $i,i+1,i+2\in \Sup(w)$, and two of these letters form a
2-cycle in $w$.
%transposition in $\Pair(w)$.
%This is an orbit of order 3.
\item[6.] $i,i+1,i+2\in \Sup(w)$, and no two of these letters form
a 2-cycle in $w$.
%a transposition in $\Pair(w)$.
%This is an orbit of order 6.
\end{itemize}

%\noindent{\bf Example.} Let $n=7$, $i=1$ and $k=2$. Then
%$(4,5)(6,7)$ is a representative of the first orbit; $(1,4)(5,6)$
%is a representative of the second; $(1,2)(4,5)$ of the third;
%$(1,4)(2,5)$ of the fourth; $(1,2)(3,4)$ of the fifth; for $k=2$
%there is no orbit of sixth type; for $k=3$, $(1,4)(2,5)(3,6)$ is
%such a representative.
%
%\medskip

Note that an orbit of the first type is of order one; orbits of
the second, third and fifth type are of order three; and orbits of
the fourth and sixth type are of order six. Moreover, by
Lemma~\ref{t.ell}, orbits of the same order form isomorphic intervals in
the weak involutive order (see Definition~\ref{involutive}). In
particular, all orbits of order six have a representative $w$ of
minimal involutive length, such that the orbit has the form :
\begin{equation}\label{e.order6}
\begin{matrix}
 & &s_is_{i+1}s_iws_is_{i+1}s_i& &  \\
 &\swarrow&                    &\searrow& \\
s_is_{i+1}ws_{i+1}s_i& & & &s_{i+1}s_iws_is_{i+1}\\
\downarrow & & & &\downarrow\\
s_{i+1}ws_{i+1}& & & & s_iws_i\\
 &\searrow&                    &\swarrow& \\
 & &w& &
\end{matrix}
\end{equation}
All orbits of order three are linear posets:
\begin{equation}\label{e.linear-1}
w <_I s_iws_i <_I s_{i+1}s_iws_is_{i+1}
\end{equation}
or
\begin{equation}\label{e.linear-2}
w <_I s_{i+1}ws_{i+1} <_I s_is_{i+1}ws_{i+1}s_i.
 \end{equation}

\smallskip

Thus %, as in the proof of Theorem~\ref{t.representation},
 the analysis %the six types may be
is reduced into three cases.

\smallskip

\noindent{\bf Case (a).} An orbit of order six.
By (\ref{e.Hecke-action}) and (\ref{e.order6}), the representation matrices of the
generators with respect to the ordered basis $C_w$, $C_{s_iws_i}$,
$C_{s_{i+1}s_iws_is_{i+1}}$, $C_{s_is_{i+1}s_iws_is_{i+1}s_i}$,
$C_{s_{i+1}ws_{i+1}}$, $C_{s_is_{i+1}ws_{i+1}s_i}$ are :
$$
\rho_q(T_i)= \left(\begin{matrix}
1-q       & 1        & 0         & 0  & 0    & 0\\
q         & 0        & 0         & 0  & 0    & 0\\
0         & 0        & 1-q       & 1  & 0    & 0  \\
0         & 0        & q         & 0  & 0    & 0  \\
0         & 0        & 0         & 0  & 1-q  & 1  \\
0         & 0        & 0         & 0  & q    & 0  \\
\end{matrix}\right)
$$
and
$$
\rho_q(T_{i+1})= \left(\begin{matrix}
1-q       & 0        & 0         & 0   & 1  & 0\\
0         & 1-q      & 1         & 0   & 0  & 0\\
0         & q        & 0         & 0   & 0  & 0  \\
0         & 0        & 0         & 0   & 0  & q  \\
q         & 0        & 0         & 0   & 0  & 0  \\
0         & 0        & 0         & 1   & 0  & 1-q \\
\end{matrix}\right)
$$
It is easy to verify that indeed
$$\rho_q(T_i)\rho_q(T_{i+1})\rho_q(T_i)=\rho_q(T_{i+1})\rho_q(T_i)\rho_q(T_{i+1}).
$$

\medskip

\noindent{\bf Case (b).} An orbit of order three. Without loss of
generality, the orbit is of type~(\ref{e.linear-1});
the analysis of type~(\ref{e.linear-2}) is analogous.
%Namely, $w <_I s_iws_i <_I s_{i+1}s_iws_is_{i+1}$.
Then $s_{i+1}ws_{i+1}=w$ and
$s_i(s_{i+1}s_iws_is_{i+1})s_i=s_{i+1}s_iws_is_{i+1}$. By
(\ref{e.iff8}), $s_{i+1}\in \Des(w)$ if and only if
$s_i\in\Des(s_{i+1}s_iws_is_{i+1})$, see second proof of
Theorem~\ref{t.representation}.
% An analogous statement holds for the case (\ref{e.linear-2}).

Given the above, by~(\ref{e.Hecke-action}), the representation
matrices of the generators with respect to the ordered basis
$w <_I s_iws_i <_I s_{i+1}s_iws_is_{i+1}$ are
$$
\rho_q(T_i)=\left(
\begin{matrix}
1-q        & 1         & 0       \\
q          & 0         & 0   \\
0          & 0         & x  \\
\end{matrix}\right)
$$
and
$$
\rho_q(T_{i+1})= \left(\begin{matrix}
x         & 0         & 0       \\
0         & 1-q       & 1   \\
0         & q         & 0   \\
\end{matrix}\right),
$$
where $x\in\{1,-q\}$. These matrices satisfy the required braid relation.

\medskip

\noindent{\bf Case (c).} An orbit of order one.
Then $s_iws_i=w$, $s_{i+1}ws_{i+1}=w$ %Furthermore, in this case
%$Sup(w)\cap\{i,i+1,i+2\}=\emptyset$, so that
and $s_i,s_{i+1}\not\in\Des(w)$. By~(\ref{e.Hecke-action}),
$\rho_q(T_i)\rho_q(T_{i+1})\rho_q(T_i)C_w=\rho_q(T_{i+1})\rho_q(T_i)\rho_q(T_{i+1})C_w=C_w$,
completing the proof of the third relation.

\medskip

The proof of the other two relations is easier and will be left to the reader.

\bigskip

\noindent{\bf Part 2: Proof of Theorem~\ref{t.Hecke-model}(2).}
%$\rho_q$ is a model.
Consider the Hecke algebra $H_n(q)$ as the algebra over $\QQ(q^{1/2})$ spanned by
$\{T_v|\ v\in S_n\}$ with the multiplication
rules
$$
T_vT_u=T_{vu} \qquad\text{if } \ell(vu)=\ell(v)+\ell(u)
$$
and
$$
(T_s+q)(T_s-1)=0 \qquad (\forall s\in S).
$$
By Lusztig's version of Tits' deformation theorem~\cite[Theorem
3.1]{L}, the group algebra of $S_n$ over $\QQ(q^{1/2})$ may be
embedded in $H_n(q)$. In particular, every element $w\in S_n$ may
be expressed as a linear combination
$$
w=\sum\limits_{v\in S_n} m_{v,w}(q^{1/2}) T_v,
$$
where $m_{v,w}$ is a rational function of $q^{1/2}$.

It follows that $\rho_q$ may be considered as an
$S_n$-representation, via
$$
\rho_q(w):=\sum\limits_{v\in S_n} m_{v,w}(q^{1/2}) \rho_q( T_v)\qquad(\forall w\in S_n).
$$
The resulting character values $\rho_q(w)$ are rational functions of $q^{1/2}$.
By discreteness of the $S_n$ character values,
%the character is constant in ``small" generic neighborhoods and thus
%constant whenever it is defined.
each such function is locally constant wherever it is defined,
and is thus constant globally.

%On the other hand, by
By Theorem~\ref{t.model}, $\rho_q|_{q=1}=\rho$ is a
model for the group algebra of $S_n$. This completes the proof.

\qed

%\medskip
%\noindent{\bf Note:} A second proof of
%Theorem~\ref{t.representation} may be obtained by substituting $q=1$
%in the proof of Part 1.

%------------------------------------------------------------------------------
\subsection{Proof of Proposition~\ref{t.Hecke-character}}
%------------------------------------------------------------------------------

Let $\SYT_n$ be the set of all standard Young tableaux of order
$n$, and let $\SYT(\lambda)\subseteq \SYT_n$ be the subset of
standard Young tableaux of shape $\lambda$. For each partition
$\lambda$ of $n$, fix a standard Young tableau $P_\lambda\in
\SYT(\lambda)$. By~\cite[Theorem 4]{Ro1}, the value of the
irreducible $H_n(q)$-character $\chi^\lambda_q$ at $T_\mu$ is
$$
\chi^\lambda_q(T_\mu) = \sum\limits_{\{w\mapsto (P_\lambda,Q)|\
w\text{ is $\mu$-unimodal and $Q\in \SYT(\lambda)$} \}}
(-q)^{\des(w)},
$$
where the sum runs over all permutations $w\in S_n$ which are
mapped under the Robinson-Schensted (RS) correspondence to
$(P_\lambda,Q)$ for some $Q\in\SYT(\lambda)$.
%where $K$ is a Knuth class of shape $\lambda$; namely, l, and let $K_P$
%consists of all permutations, which are mapped under RSK
%correspondence, to a pair of tableaux $(P,Q)$.
By~\cite[Lemma 7.23.1]{EC2}, the descent set of $w\in S_n$, which
is mapped under RS to $(P_\lambda,Q)$, is determined by $Q$.
Hence
$$
\Tr\ \rho_q(T_\mu)= \sum\limits_\lambda \chi^\lambda_q(T_\mu))=
\sum\limits_\lambda \sum\limits_{\{w\mapsto (P_\lambda,Q)|\
w\text{ is $\mu$-unimodal
 and $Q\in \SYT(\lambda)$} \}} (-q)^{\des(w)}
$$
$$
= \sum\limits_\lambda \sum\limits_{\{w\mapsto (Q,Q)|\ w\text{ is
$\mu$-unimodal  and $Q\in \SYT(\lambda)$} \}} (-q)^{\des(w)}
$$
$$
= \sum\limits_{\{w\mapsto (Q,Q)|\ \text{ $Q\in \SYT_n$ and $w$ is
$\mu$-unimodal} \}} (-q)^{\des(w)}=\sum\limits_{\{w\in I_n|\
w\text{ is $\mu$-unimodal}\}} (-q)^{\des(w)}.
$$
The last equality follows from the well known property of the RS
correspondence: $w\mapsto (P,Q)$ if and only if $w^{-1}\mapsto
(Q,P)$~\cite[Theorem 7.13.1]{EC2}. Thus $w$ is an involution if
and only if $w\mapsto (Q,Q)$ for some $Q\in\SYT_n$.

\qed

%------------------------------------------------------------------------------
\section{Remarks and Questions}
%------------------------------------------------------------------------------

%------------------------------------------------------------------------------
\subsection{Classical Weyl Groups}
%------------------------------------------------------------------------------

Let $B_n$ be the Weyl group of type $B$,
$S^B$ its set of simple reflections,
$I_n^B$ its set of involutions,
and $V_n^B:={\rm{span}}_\QQ\{C_w\ |\ w\in I_n^B\}$ a vector space over $\QQ$
formally spanned by the involutions. Recall that
$B_n=\ZZ_2\wr S_n$, so that each element $w\in B_n$ is identified
with a pair $(v,\sigma)$, where $v\in \ZZ_2^n$ and $\sigma\in S_n$.
Denote $|w|:=\sigma$.

Define a map $\rho^B:S^B\to GL(V_n)$ by
%\begin{equation}\label{B-signed conjugation}
$$
\rho^B(s)C_w := {\rm{sign}}(s;\ w)\cdot C_{sws}\qquad (\forall s\in
S^B, w\in I_n^B)
$$
%\end{equation}
where, for $s= s_0=((1,0,\dots,0), id)$, the exceptional Coxeter generator, the sign is
%defined as in (\ref{e.sign-definition}); namely,
%\begin{equation}\label{B1-sign-definition2}
$$
{\rm{sign}}(s_0;\ w):=
\begin{cases}
   -1,   & \text{ if } sws=w \text { and } s_0\in \Des(w); \\
   1, & \text{ otherwise, }
\end{cases}
$$
%\end{equation}
and for a generator $s\ne s_0$ the sign is
%\begin{equation}\label{B1-sign-definition}
$$
{\rm{sign}}(s;\ w):=
\begin{cases}
   -1,   & \text{ if } sws=w \text { and } s\in \Des(|w|); \\
   1, & \text{ otherwise. }
\end{cases}
$$
%\end{equation}
%Here $|w|=(|w(1)|,\dots,|w(n)|)\in S_n$, the permutation obtained
%by ignoring the signs.

\begin{theorem}%\label{B-representation}
$\rho^B$ is a Gelfand model for $B_n$.
\end{theorem}

\noindent %A proof will be given elsewhere.
A proof is given in~\cite{APR}.

\medskip

Models for classical Weyl groups of type $D_n$ for odd $n$ were
constructed in~\cite{B, A3}. These constructions fail for even
$n$. %One may guess that a signed conjugation
%%,where the sign is
%%defined as in (\ref{e.sign-definition}),
%will give a model for
%$D_{2n+1}$.
%\begin{question}
%Is there a signed conjugation (or a representation of type
%$\rho_sC_w=a_{s,w} C_w + b_{s,w} C_{sws}$) which gives a model for
%$D_{2n}$ ?
%\end{question}
A natural question is whether there exists a signed conjugation
(or a representation of type $\rho_sC_w=a_{s,w} C_w + b_{s,w}
C_{sws}$) which gives a model for $D_{2n}$.
% Let $S^D$ - the set of
%simple reflections in $D_n$, $I_n^D$ - the set of involutions in
%$D_n$ and $V_n^D:={\rm{span}}_\QQ\{C_w\ |\ w\in I_n^D\}$ be a
%vector space over $\QQ$ spanned by the involutions. Define a map
%$\rho^D:S^D\to GL(V_n)$ by
%\begin{equation}\label{D-signed conjugation}
%\rho(s)C_w = {\rm{sign}}(s;\ w)\cdot C_{sws}\qquad (\forall s\in
%S^D, w\in I_n^D)
%\end{equation}
%where the sign is defined as in (\ref{e.sign-definition}). Then
%
%\begin{question}\label{D-representation} Is $\rho^D$  a
%$D_n$-model ?
%\end{question}
It is also desired to find representation matrices for the models
of the  Hecke algebras of types $B$ and $D$ which specialize at
$q=1$ to models of the corresponding group algebra.
%
%\begin{question} Find representation matrices for the
%models of the  Hecke algebras of types $B$ and $D$ which
%specialize at $q=1$ to models of the corresponding group algebra.
%\end{question}

\medskip

We conclude with the following
%interesting
questions regarding an arbitrary Coxeter group $W$.

\begin{question}
Find a signed conjugation
%(or a representation of type
%$\rho_sC_w=a_{s,w} C_w + b_{s,w} C_{sws}$),
which gives a Gelfand model for $W$;
Find %a signed conjugation (or
a representation of the form $\rho_sC_w=a_{s,w} C_w + b_{s,w}
C_{sws}$, which gives a Gelfand model for the
%(finite)
 Hecke algebra of $W$.
\end{question}

\begin{question}
Find a character formula for the Gelfand model of the  Hecke algebra of $W$.
\end{question}

%------------------------------------------------------------------------------
%\section{A Second Proof of Theorem~\ref{t.model}}
\section{Appendix}
%------------------------------------------------------------------------------

This section was added in proof.

\bigskip

First, it should be acknowledged that an equivalent reformulation
of Theorem~\ref{t.model}, with a different proof, was given
by Kodiyalam and Verma~\cite{Verma}.

\bigskip

A third proof of Theorem~\ref{t.model}, along the lines of~\cite{Saxl},
was suggested by an anonymous referee. Here is a brief outline.

\medskip

%The starting point is a well known result~\cite[Ch. I \S 8 Ex 6,
%and Ch. VII (2.4)]{Md}.

Let $\chi^{\emptyset,(n)}$ denote the one dimensional character of
$B_n$ given by the parity of the negative signs, and consider the
natural embedding of $B_n=\ZZ_2\wr S_n$ into $S_{2n}$. Then
$$
\chi^{\emptyset,(n)}\uparrow^{S_{2n}}_{B_n}=
\sum\limits_{\lambda\vdash n}\chi^{(2\cdot \lambda)'},
$$
where the sum on the right hand side runs through all partitions of $2n$
with even columns only.  See, for example, ~\cite[Ch.~I \S 8 Ex.~6,
and Ch.~VII (2.4)]{Md}. Combining this with the Littlewood-Richardson rule
implies that
$$
((\chi^{\emptyset,(k)}\uparrow^{S_{2k}}_{B_k})\otimes
1_{S_{n-2k}})\uparrow_{S_{2k}\times S_{n-2k}}^{S_n}
$$
is a multiplicity free sum of all irreducible Specht modules
indexed by partitions with exactly $n-2k$ odd columns.

A natural basis for this representation is given by involutions
with $n-2k$ fixed points. Finally, it is straightforward to show
that the action of a Coxeter generator $s_i$ on this basis is
identical with the signed conjugation defined in (\ref{e.sign-definition}).

\qed

\medskip

\begin{corollary}
The signed conjugation $\rho$, when restricted to the conjugacy
class of involutions with $n-2k$ fixed points, is a multiplicity
free sum of all irreducible Specht modules indexed by partitions
with exactly $n-2k$ odd columns.
\end{corollary}

This gives an algebraic proof to the following combinatorial
result: The number of involutions with $n-2k$ fixed points is
equal to the number of standard Young tableaux of shapes with
exactly $n-2k$ odd columns.

Another proof for this enumerative fact may be obtained using the
Robinson-Schensted (RS) correspondence. One concludes that  the
restriction of the signed conjugation $\rho$ to conjugacy classes
of involutions
%is compatible with the RS correspondence. In other
%words, the decomposition of $\rho$ into irreducibles
%is essentially given
is compatible with the RS correspondence; namely, $S^\lambda$ is
a factor of the restriction of $\rho$ to the conjugacy class of
cycle type $2^k 1^{n-2k}$ if and only if $\lambda$ is the shape of
some pair of (equal) standard Young tableaux corresponding to an involution of
this cycle type.

%the irreducible Specht modules which appear in the decomposition
%are exactly those indexed by the partitions obtained via the
%Robinson-Schensted correspondence.

%An explicit description of the irreducible modules which occur
%when the signed conjugation is restricted to conjugacy classes of
%involutions.

%\begin{corollary}
%The signed conjugation $\rho$, when restricted to the conjugacy
%classes of involution with $n-2k$ fixed points, decomposes into
%the irreducible Specht modules indexed by the partitions obtained
%via the Robinson-Schensted correspondence from these involutions.
%\end{corollary}

%------------------------------------------------------------------------------
\section{Acknowledgements}
The authors thank Arkady Berenstein,  Steve Shnider and Richard
Stanley for stimulating discussions and references. We also thank
an anonymous referee for his contributions and suggestions.
%------------------------------------------------------------------------------

%------------------------------------------------------------------------------

\end{document}